\documentclass{amsart}

\DeclareMathSymbol{\twoheadrightarrow}
{\mathrel}{AMSa}{"10}

\def\bchi{{\mathbf \chi}}
\def\bphi{{\mathbf \phi}}
\def\Q{{\mathbf Q}}
\def\Z{{\mathbf Z}}

\def\R{{\mathbf R}}
\def\RR{{\mathfrak R}}
\def\F{{\mathbf F}}

\def\A{{\mathbf A}}
\def\Sn{{\mathbf S}_n}
\def\An{{\mathbf A}_n}

\def\Gal{\mathrm{Gal}}
\def\Perm{\mathrm{Perm}}

\def\End{\mathrm{End}}
\def\Aut{\mathrm{Aut}}

\def\I{\mathrm{Id}}

\def\fchar{\mathrm{char}}

\def\GL{\mathrm{GL}}
\def\Sz{\mathrm{Sz}}
\def\L{\mathrm{L}}
\def\U{\mathrm{U}}
\def\B{\mathrm{B}}
\def\PGL{\mathrm{PGL}}
\def\SL{\mathrm{SL}}
\def\PSL{\mathrm{PSL}}
\def\PSU{\mathrm{PSU}}
\def\Sp{\mathrm{Sp}}

\def\dim{\mathrm{dim}}

\def\B{{B}}

\newtheorem{thm}{Theorem}[section]
\newtheorem{lem}[thm]{Lemma}

\theoremstyle{definition}
\newtheorem{defn}[thm]{Definition}
\newtheorem{ex}[thm]{Example}

\newtheorem{rem}[thm]{Remark}
\newtheorem{rems}[thm]{Remarks}
\title[Hyperelliptic jacobians and Steinberg representations]
{Hyperelliptic jacobians without complex multiplication and
Steinberg representations in positive characteristic}
\author[Yuri G. Zarhin]{Yuri G. Zarhin}
\address{Department of Mathematics, Pennsylvania State University,
University Park, PA 16802, USA} \email{zarhin\char`\@math.psu.edu}

\begin{document}
\begin{abstract}
In his previous papers ~\cite{ZarhinMRL,ZarhinMRL2,ZarhinMMJ} the
author proved that in characteristic $\ne 2$ the jacobian $J(C)$
of a hyperelliptic curve $C: y^2=f(x)$ has only trivial
endomorphisms over an algebraic closure $K_a$ of the ground field
$K$ if the Galois group $\Gal(f)$ of the irreducible polynomial
$f(x) \in K[x]$ is either the symmetric group $\Sn$ or the
alternating group $\A_n$. Here $n\ge 9$ is the degree of $f$. The
goal of this paper is to extend this result to the case of certain
``smaller'' doubly transitive Galois groups. Namely, we treat the
infinite series $n=2^m+1, \Gal(f)=\L_2(2^m):=\PSL_2(\F_{2^m})$,
$n=2^{4m+2}+1, \Gal(f)=\Sz(2^{2m+1})= {^2\B_2}(2^{2m+1})$ and
$n=2^{3m}+1, \Gal(f)=\U_3(2^m):=\PSU_3(\F_{2^m})$.



\end{abstract}
\maketitle
\section{Introduction}
Let $K$ be a field and $K_a$ its algebraic closure. Assuming that
$\fchar(K)=0$, the author \cite{ZarhinMRL} proved that  the
jacobian $J(C)=J(C_f)$  of a hyperelliptic curve
$$C=C_f:y^2=f(x)$$
has only trivial endomorphisms over  $K_a$ if the Galois group
$\Gal(f)$ of the irreducible polynomial $f \in K[x]$ is ``very
big". Namely, if $n=\deg(f) \ge 5$ and $\Gal(f)$ is either the
symmetric group $\Sn$ or the alternating group $\An$
 then the ring $\End(J(C_f))$ of $K_a$-endomorphisms of $J(C_f)$ coincides with $\Z$.
 Later the author ~\cite{ZarhinTexel,ZarhinPAMS} proved that $\End(J(C_f))=\Z$ for  infinite series
 $\Gal(f)=\L_2(2^m):=\PSL_2(\F_{2^m})$ and $n=2^{m}+1$ (with
$\dim(J(C_f))=2^{m-1}$), $\Gal(f)=$
 the Suzuki group $\Sz(2^{2m+1})={^2\B_2}(2^{2m+1})$ and $n=2^{2(2m+1)}+1$
(with $\dim(J(C_f))=2^{4m+1}$),
$\Gal(f)=\U_3(2^m):=\PSU_3(\F_{2^m})$ and $n=2^{3m}+1$ (with
$\dim(J(C_f))=2^{3m-1}$). We refer the reader to
~\cite{Mori1,Mori2,Katz1,Katz2,Masser,KS,ZarhinMRL,ZarhinMRL2,ZarhinMMJ,ZarhinP}
for a discussion of known results about, and examples of,
hyperelliptic jacobians without complex multiplication.

When $\fchar(K)>2$, the author \cite{ZarhinMMJ} proved  that
$\End(J(C_f))=\Z$ if $n \ge 9$ and $\Gal(f)=\Sn$ or $\An$. The aim
of the present paper is to extend this result to the case of
already mentioned series of doubly transitive Galois groups
$\L_2(2^m),\Sz(2^{2m+1})$ and $\U_3(2^m)$. Notice that it is known
~\cite{ZarhinTexel,ZarhinPAMS} that in this case  either
$\End(J(C))=\Z$ or $J(C)$ is a supersingular abelian variety and
the real problem is how to prove that $J(C)$ is {\sl not}
supersingular.

\begin{rem}
The groups $\L_2(2^m),\U_3(2^m)$ (with $m \ge 2$) and $\Sz(2^{2m+1})$
 constitute an interesting important class of finite simple groups called simple
 Bender groups ~\cite[Chapter 2]{GLS}. Our main theorem deals with all these groups
 except $\L_2(4)\cong {\mathrm A_5}$.
\end{rem}

\section{Main result}
\label{mainr}
Throughout this paper we assume that $K$ is a field of prime characteristic
$p$
different from $2$. We fix its algebraic closure $K_a$
and write $\Gal(K)$ for the absolute Galois group $\Aut(K_a/K)$.

\begin{thm}
\label{main}
Let $K$ be a field with $p=\fchar(K) > 2$,
 $K_a$ its algebraic closure,
$f(x) \in K[x]$ an irreducible separable polynomial of  degree $n
\ge 9$. Let us assume  that $n$ and the Galois group $\Gal(f)$ of
$f$ enjoy one of the following properties:
\begin{enumerate}
\item[(i)]
There exists a positive integer $m\ge 3$ such that
 $n=2^m+1$ and  $\Gal(f)$ contains
 a subgroup isomorphic to $\L_2(2^m)$;
\item[(ii)]
There exists a  positive integer $m\ge 1$ such that
$n=2^{2(2m+1)}+1$ and $\Gal(f)$  contains
 a subgroup isomorphic to $\Sz(2^{2m+1})$;
\item[(iii)]
There exists a positive integer $m\ge 2$ such that $n=2^{3m}+1$
and $\Gal(f)$ contains
 a subgroup isomorphic  to $\U_3(2^m)$.
\end{enumerate}

 Let $C_f$ be the hyperelliptic curve
$y^2=f(x)$. Let $J(C_f)$ be its jacobian, $\End(J(C_f))$ the ring
of $K_a$-endomorphisms of $J(C_f)$. Then $\End(J(C_f))=\Z$.
\end{thm}

\begin{rem}
\label{redA}
 Replacing $K$ by a suitable finite separable
extension, we may assume in the course of the proof of Theorem
\ref{main} that $\Gal(f)=\L_2(2^m),\Sz(2^{2m+1})$ or $\U_3(2^m)$
respectively. Taking into account that all these groups are simple
non-abelian and replacing $K$ by its abelian extension obtained by
adjoining to $K$ all $2$-power roots of unity, we may also assume
that $K$ contains all $2$-power roots of unity.
\end{rem}

As was already pointed out, in light of results of
~\cite{ZarhinTexel,ZarhinPAMS} and Remark \ref{redA}, our
 Theorem \ref{main} is an immediate corollary of the following auxiliary
statement.

\begin{thm}
\label{main2} Suppose   $n$ is an odd integer.  Suppose $K$ is a
field, $\fchar(K)=p \ne 2$ and $K$ contains all $2$-power roots of
unity. Suppose that $f(x) \in K[x]$ is a separable polynomial of
degree $n$, whose Galois group $\Gal(f)$ enjoys one of the
following properties:
\begin{enumerate}
\item[(i)]
There exists a positive integer $m\ge 3$ such that
 $n=2^m+1$ and  $\Gal(f)=\L_2(2^m)$;
\item[(ii)]
There exists a  positive integer $m$ such that
$n=2^{2(2m+1)}+1$ and $\Gal(f)=\Sz(2^{2m+1})$;
\item[(iii)]
There exists a positive integer $m\ge 2$ such that $n=2^{3m}+1$
and $\Gal(f)=\U_3(2^m)$.
\end{enumerate}

 Let $C$ be
the hyperelliptic curve $y^2=f(x)$ of genus $g=\frac{n-1}{2}$ over
$K$ and let $J(C)$  be the jacobian of $C$.

Then $J(C)$ is not a supersingular abelian variety.
\end{thm}

\begin{ex}
Let $k$ be an algebraically closed field of characteristic $7$.
Let $K=k(z)$ be the field of rational functions in variable $z$
with constant field $k$. We write $\overline{k(z)}$ for an
algebraic closure of $k(z)$. According to Abhyankar \cite{Ab7},
the Galois group of the polynomial
$$x^9-z x^7+1 \in k(z)[x]=K[x]$$
is $\L_2(8)$ (see also ~\cite[\S 3.3, Remarque 2(a)]{SerreB}).
Hence the four-dimensional jacobian of the hyperelliptic curve
$y^2=x^9-z x^7+1$ has no nontrivial endomorphisms over
$\overline{k(z)}$.
\end{ex}

We prove Theorem \ref{main2} in  Section \ref{pmain2}.

\section{Permutation groups, permutation modules and very simple representations}
\label{permute}

Let $B$ be a finite set consisting of $n \ge 5$ elements. We write
$\Perm(B)$ for the group of permutations of $B$. A choice of
ordering on $B$ gives rise to an isomorphism $\Perm(B) \cong \Sn$.
Let $G$ be a  subgroup of $\Perm(B)$. For each $b \in B$ we write
$G_b$ for the stabilizer of $b$ in $G$; it is a subgroup of $G$.

\begin{rem}
\label{transitive} Assume that the action of $G$ on $B$ is
transitive. It is well-known that each $G_b$ is  of index $n$ in
$G$ and
 all the $G_b$'s are conjugate in $G$.
 Each conjugate of $G_b$ in $G$ is the stabilizer of a point of $B$.
 In addition, one may identify the $G$-set $B$ with the set of cosets $G/G_b$
with the standard action by $G$.
\end{rem}

Let $\F$ be a field. We write $\F^B$ for the $n$-dimensional
$\F$-vector space of maps $h:B \to \F$. The space $\F^B$ is
provided with a natural action of $\Perm(B)$ defined as follows.
Each $s \in \Perm(B)$ sends a map
 $h:B\to \F$ into  $sh:b \mapsto h(s^{-1}(b))$. The permutation module $\F^B$
contains the $\Perm(B)$-stable hyperplane
$$(\F^B)^0=
\{h:B\to\F\mid\sum_{b\in B}h(b)=0\}$$ and the $\Perm(B)$-invariant
line $\F \cdot 1_B$ where $1_B$ is the constant function $1$.

Clearly, $(\F^B)^0$ contains $\F \cdot 1_B$ if and only if
$\fchar(\F)$ divides $n$. If this is {\sl not} the case then there
is a $\Perm(B)$-invariant splitting $\F^B=(\F^B)^0 \oplus \F \cdot
1_B$.

Clearly, $\F^B$ and $(\F^B)^0$  carry natural structures of
$G$-modules. Their characters depend only on the characteristic of
$\F$.

Let us consider the case of $\F=\Q$. Then the character of $\Q^B$
 is called the {\sl permutation character} of $B$. Let us denote by
$\bchi=\bchi_B:G \to \Q$ the character of $(\Q^B)^0$.
 Clearly, $1+\bchi$ is the permutation character of $B$.
It is also clear that the representation of $G$ in $(\Q^B)^0$ is
{\sl orthogonal}.

Now, let us consider the case of $\F=\F_2$. If  $n$ is even then
let us define the $\Perm(B)$-module $Q_B:=(\F_2^B)^0/(\F_2 \cdot
1_B)$. If $n$ is odd then let us put $Q_B:=(\F_2^B)^0$.

\begin{rem}
\label{lift}
 Clearly, $Q_B$ is a faithful $G$-module (recall that
$n \ge 5$). If $n$ is odd then $\dim_{\F_2}(Q_B)=n-1$ and one may
view the $\Q[G]$-module $(\Q^B)^0$ as a {\sl lifting} to
characteristic zero of the $\F_2[G]$-module $(\Q^B)^0$ . If $n$ is
even then $\dim_{\F_2}(Q_B)=n-2$.
\end{rem}

 Let $G^{(2)}$ be the set of $2$-regular elements of $G$.
 Clearly, the Brauer character of the $G$-module $\F_{2}^B$
 coincides with the restriction of $1+\bchi_B$ to $G^{(2)}$.
 This implies easily that the Brauer character of the $G$-module $(\F_{2}^B)^0$
 coincides with the restriction of   $\bchi_B$ to $G^{(2)}$.

\begin{rem}
\label{Bcharacter}
 Let us denote by
$\bphi_B=\bphi$
 the Brauer character of the $G$-module $Q_B$.
 One may easily check that $\bphi_B$ coincides with the restriction of
 $\bchi_B$ to $G^{(2)}$ if  $n$ is odd and with the restriction of
$\bchi_B-1$ to $G^{(2)}$ if  $n$ is even.
\end{rem}

We refer to ~\cite{ZarhinTexel,ZarhinMRL2,ZarhinMMJ} for a
discussion of the following definition.

\begin{defn}
Let $V$ be a vector space over a field $\F$, let $G$ be a group
and $\rho: G \to \Aut_{\F}(V)$ a linear representation of $G$ in
$V$. We say that the $G$-module $V$ is {\sl very simple} if it
enjoys the following property:

If $R \subset \End_{\F}(V)$ is an $\F$-subalgebra containing the
identity operator $\I$ such that

 $$\rho(\sigma) R \rho(\sigma)^{-1} \subset R \quad \forall \sigma \in G$$
 then either $R=\F\cdot \I$ or $R=\End_{\F}(V)$.
\end{defn}

\begin{rems}
\label{image}
\begin{enumerate}
\item[(i)]
If $G'$ is a subgroup of $G$ and the $G'$-module $V$ is very
simple then obviously the $G$-module $V$ is also very simple.

\item[(ii)]
A  very simple module is absolutely simple (see \cite[Remark
2.2(ii)]{ZarhinTexel}).

\item[(iii)]
Clearly, the $G$-module $V$ is very simple if and only if the
corresponding $\rho(G)$-module $V$ is very simple. This implies
easily that if $H \twoheadrightarrow G$ is a surjective group
homomorphism then the $G$-module $V$ is  very simple if and only
if the corresponding $H$-module $V$ is  very simple.

\item[(iv)]
Let $G'$ be a normal subgroup of $G$. If $V$ is a very simple
$G$-module then either $\rho(G') \subset \Aut_{k}(V)$ consists of
scalars (i.e., lies in $k\cdot\I$) or the $G'$-module $V$ is
absolutely simple. See ~\cite[Remark 5.2(iv)]{ZarhinMMJ}.

\item[(v)]
Suppose $F$ is a discrete valuation field with  valuation ring
$O_F$, maximal ideal $m_F$ and residue field $k=O_F/m_F$. Suppose
$V_F$ a finite-dimensional $F$-vector space, $\rho_F: G \to
\Aut_F(V_F)$ a $F$-linear representation of $G$. Suppose $T$ is a
$G$-stable $O_F$-lattice in $V_F$ and the corresponding
$k[G]$-module $T/m_F T$ is isomorphic to $V$. Assume that the
$G$-module $V$ is very simple. Then the $G$-module $V_F$ is also
very simple. See ~\cite[Remark 5.2(v)]{ZarhinMMJ}.
\end{enumerate}
\end{rems}

\begin{thm}
\label{very}
 Suppose one of the following conditions holds:
\begin{enumerate}
\item[(i)]
There exists a positive integer $m\ge 3$ such that
 $n=2^m+1$ and  $G \cong\L_2(2^m)$;
\item[(ii)]
There exists a  positive integer $m$ such that $n=2^{2(2m+1)}+1$
and $G \cong \Sz(2^{2m+1})$;
\item[(iii)]
There exists a positive integer $m\ge 2$ such that $n=2^{3m}+1$
and $G\cong\U_3(2^m)$.
\end{enumerate}
Then the $G$-module $Q_B$ is very simple.
\end{thm}

\begin{proof}
See Th. 7.10 and Th. 7.11 of \cite{ZarhinTexel} and Th. 4.4 of
\cite{ZarhinPAMS}.
\end{proof}

The following statement provides a criterion for a representation
over $\F_2$ to be  very simple.

\begin{thm}
\label{Very3} Suppose that a positive integer $N>1$ and a group
$H$ enjoy the following properties:
\begin{itemize}
\item
$H$ does not contain a subgroup of index dividing $N$ except $H$
itself.

\item
Let $N=ab$ be a factorization of $N$ into a product of two
positive integers $a>1$ and $b>1$. Then either there does not
exist an absolutely simple $\F_2[H]$-module of $\F_2$-dimension
$a$ or there does not exist an absolutely simple $\F_2[H]$-module
of $\F_2$-dimension $b$.
\end{itemize}

Then each absolutely simple $\F_2[H]$-module of $\F_2$-dimension
$N$ is very simple.
\end{thm}

\begin{proof}
This is Corollary 4.12 of \cite{ZarhinTexel}.
\end{proof}

\begin{thm}
\label{Lnq} Suppose that there exist a positive integer $m>2$ and
an odd power prime $q$ such that $n=\frac{q^m-1}{q-1}$. Suppose
$G$ is a subgroup of $\Sn$ and contains  a subgroup isomorphic to
$\L_m(q)$. Then the $G$-module $Q_B$ is very simple.
\end{thm}

\begin{proof}

 In light of Remark \ref{image}(i), we may assume that $G=\L_m(q)$.

Assume that $(m,q) \ne (4,3)$. Then the assertion of Theorem
\ref{Lnq} follows easily from ~\cite[Remark 4.4 and Corollary 5.4
]{ZarhinMMJ}. (The proof in \cite{ZarhinMMJ} was based on a result
of Guralnick \cite{GurTiep}.)

So, we may assume that
 $m=4,q=3$. We have $n=\#(B)=40$ and $\dim_{\F_2}(Q_B)=38$.
According to the Atlas ~\cite[pp. 68-69]{Atlas}, $G=\L_4(3)$ has
two conjugacy classes of maximal subgroups of index $40$. All
other maximal subgroups have index greater than $40$.
 Therefore all subgroups of $G$ (except $G$ itself) have index greater than $39>38$.
This implies that each action of $G$ on $B$ is transitive. The
corresponding permutation character (in notations of \cite{Atlas})
coincides (in both cases) with $1+\chi_4$, i.e., $\bchi=\chi_4$.
Since $40$ is even, we need to consider the restriction of
$\bchi-1$ to the set of $2$-regular elements of $G$ and this
restriction coincides with the absolutely irreducible Brauer
character
 $\phi_4$ (in notations of \cite{AtlasB}, p. 165).
 In particular, the corresponding $G$-module $Q_B$ is absolutely simple.
It follows from the Table on p. 165 of \cite{AtlasB} that all
absolutely irreducible representations
 of $G$ in characteristic $2$ have dimension which is {\sl not}
 a strict divisor of $38$.
Combining this observation with the absence of subgroups in $G$ of
index $\le 38$, we conclude, thanks to Theorem \ref{Very3}, that
$Q_B$ is very simple. This ends the proof of Theorem \ref{Lnq}.
\end{proof}

\section{Proof of Theorem \ref{main2}}
\label{pmain2}
 So, we assume that $K$ contains all $2$-power roots
of unity, $f(x) \in K[x]$ is an irreducible separable polynomial
of degree $n=2g+1$ or $2g+2$ and $n \ge 5$. Therefore the jacobian
$J(C)$ of the hyperelliptic curve $C:y^2=f(x)$ is a
$g$-dimensional abelian variety defined over $K$. The group
$J(C)_2$ of its points of order $2$ is a $2g$-dimensional
$\F_2$-vector space provided with the natural action of $\Gal(K)$.
It is well-known (see for instance ~\cite[Sect. 5]{ZarhinTexel})
that the image of $\Gal(K)$ in $\Aut(J(C)_2)$ is canonically
isomorphic to $\Gal(f)$. Let us recall a well-known explicit
description of the $\Gal(f)$-module $J(C)_2$ (see, for instance,
\cite{Mori2} or \cite{ZarhinTexel}). Let $\RR\subset K_a$ be the
$n$-element set of roots of $f$. We view $G=\Gal(f)$ as a certain
subgroup of the group $\Perm(\RR) \cong \Sn$ of all permutations
of $\RR$. If we put $B=\RR$ then the $2g$-dimensional
$\F_2$-vector space $Q_{\RR}=(\F_2^{\RR})^0$ carries the natural
structure of faithful $G$-module. It is well-known that the
natural homomorphism $\Gal(K) \to \Aut_{\F_2}(J(C)_2)$ factors
through the canonical surjection $\Gal(K) \twoheadrightarrow
\Gal(f)=G$ and the $G$-modules $J(C)_2$ and $(\F_2^{\R})^0$ are
isomorphic ~\cite{Mori1,ZarhinTexel}.

We deduce Theorem  \ref{main2} from the following assertion.

\begin{lem}
\label{supernot} Let $F$ be a field, whose characteristic is not
$2$ and assume that $F$ contains all $2$-power roots of unity. Let
$F_a$ be an algebraic closure of $F$. Let $g$ be a positive
integer and $G$ be a finite simple non-abelian group enjoying the
following properties:
\begin{enumerate}
\item[(a)]
Let $2^r$ be the largest power of $2$ that divides the order of
$G$. Then either $2^r$ divides $2g$  or $G=\L_4(3),g=19$.
\item[(b)]
Either the Schur multiplier of $G$ is a group of odd order or
$G=\Sz(8),g=32$.
\end{enumerate}

Suppose  that $X$ is a $g$-dimensional abelian variety over $F$
such that the image of $\Gal(F)$ in $\Aut(X_2)$ is isomorphic to
$G$ and the corresponding faithful representation
$$\rho: G \hookrightarrow \Aut(X_2) \cong \GL(2g,\F_2)$$
enjoys the following properties:
\begin{enumerate}
\item[(c)]
The representation $\rho$ is {\sl very simple}; in particular, it
is absolutely irreducible;
\item[(d)]
If $(G,g) \ne (\L_4(3),19)$ then one may lift $\rho$ to an
orthogonal representation of $G$ in characteristic zero.
\end{enumerate}
Then the ring $\End(X)$ of all $F_a$-endomorphisms of $X$
coincides with $\Z$.
\end{lem}

Lemma \ref{supernot} will be proven in the next Section.

\begin{proof}[Proof of Theorem \ref{main2}]
Clearly, $n\ge 9$ is odd. Let us put
$$G=\Gal(f), g=\frac{n-1}{2}, F=K, X=J(C).$$
Notice that $n-1$ is a power of $2$. We already observed that the
image of $\Gal(K)$ in $\Aut(J(C)_2)$ is isomorphic to $\Gal(f)$.
Let us recall a well-known explicit description of the
$\Gal(f)$-module $J(C)_2$ (see, for instance, \cite{Mori2} or
\cite{ZarhinTexel}). Let $\RR\subset K_a$ be the $n$-element set
of roots of $f$. We view $G=\Gal(f)$ as a certain subgroup of the
group $\Perm(\RR) \cong \Sn$ of all permutations of $\RR$. It is
known (see for instance  ~\cite{ZarhinTexel,ZarhinPAMS}) that
under the assumptions of Theorem \ref{main2} $G$ is a doubly
transitive permutation group. If we put $B=\RR$ then the
$(n-1)$-dimensional $\F_2$-vector space $Q_{\RR}=(\F_2^{\RR})^0$
carries the natural structure of faithful $G$-module.

 It is well-known that the natural homomorphism $\Gal(K) \to \Aut_{\F_2}(J(C)_2)$
factors through the canonical surjection $\Gal(K)
\twoheadrightarrow \Gal(f)=G$ and the $G$-modules $J(C)_2$ and
$(\F_2^{\R})^0$ are isomorphic ~\cite{Mori1,ZarhinTexel}.

Now consider the case of $\F=\Q$. Recall that the $G$-module
$(\Q^{\RR})^0$ is orthogonal.

 Clearly, the $\Q[G]$-module $(\Q^{\RR})^0$ is
a lifting to characteristic zero of the $\F_2[G]$-module
$Q_{\RR}$. This proves that the condition (d) of Lemma
\ref{supernot} holds true.

On the other hand, it follows from Theorem \ref{very} that  the
$G=\Gal(f)$-module $(\F_2^{\RR})^0=J(C)_2$ is very simple. (It is
actually a {\sl Steinberg} representation of $G$ which explains
the title of this paper.) This proves that the condition (c) of
Lemma \ref{supernot} holds true. The validity of the condition
\ref{supernot}(a) follows from the known formulas for the orders
of the simple groups involved ~\cite[p. 8]{GLS1}. (In fact,
$2g=n-1$ does coincide with the largest power of $2$ dividing the
order of $G$.) It follows from the Table in ~\cite[\S 4.15A]{G}
that the order of the Schur multiplier of $G$ is an odd number
except the case $G=\Sz(32), n=65, g=32$.
  So, we may apply Lemma \ref{supernot}
to $X=J(C)$ and conclude that $J(C)$ is {\sl not} supersingular.
\end{proof}

\section{Abelian varieties without complex multiplication}
\label{main2p} We keep all the notations and assumptions of Lemma
\ref{supernot}.

We write $T_2(X)$ for the $2$-adic Tate
module of $X$ and
$$\rho_{2,X}:\Gal(F) \to \Aut_{\Z_2}(T_2(X))$$
for the corresponding $2$-adic representation. It is well-known
that $T_2(X)$ is a free $\Z_2$-module of rank $2\dim(X)=2g$ and
$X_2=T_2(X)/2 T_2(X)$ (the equality of Galois modules). Let us put
$$H=\rho_{2,X}(\Gal(F)) \subset \Aut_{\Z_2}(T_2(X)).$$
Clearly, the natural homomorphism $\bar{\rho}_{2,X}:\Gal(F) \to
\Aut(X_2)$ defining the Galois action on the points of order $2$
is the composition of $\rho_{2,X}$ and (surjective) reduction map
modulo $2$
$$\Aut_{\Z_2}(T_2(X)) \twoheadrightarrow \Aut(X_2).$$
This gives us a natural (continuous) {\sl surjection}
$$\pi:H \twoheadrightarrow \bar{\rho}_{2,X}(\Gal(F)) \cong G,$$
 whose kernel consists of elements of $1+2 \End_{\Z_2}(T_2(X))$.
 By \ref{supernot}(c), the
$G$-module $X_2$ is very simple. By  Remark \ref{image}(iii), the
$H$-module $X_2$ is also very simple. Here the structure of
$H$-module is defined on $X_2$ via
$$H\subset\Aut_{\Z_2}(T_2(X)) \twoheadrightarrow \Aut(X_2).$$
Clearly, the $\Q_2[H]$-module $V_2(X)$ is a lifting to
characteristic zero of the very simple $\F_2[H]$-module $X_2$ and
therefore is also very simple, thanks to Remark \ref{image}(v).
Here $V_2(X)=T_2(X)\otimes_{\Z_2}\Q_2$ is the $\Q_2$-Tate module
of $X$ containing $T_2(X)$ as a $H$-stable $\Z_2$-lattice. In
particular, the $\Q_2[H]$-module $V_2(X)$ is absolutely simple.

The choice of polarization on $X$ gives rise to a non-degenerate alternating
bilinear form (Riemann form) \cite{MumfordAV}
$$e: V_{2}(X) \times V_2(X) \to \Q_2(1) \cong \Q_2.$$
Since $F$ contains all $2$-power roots of unity, $e$ is $\Gal(F)$-invariant
and therefore is $H$-invariant. In particular,
$$H \subset \Sp(V_2(X),e)\subset\SL(V_2(X)).$$
Here $\Sp(V_2(X),e)$ is the symplectic group attached to $e$. In
particular, the $H$-module $V_2(X)$ is symplectic.

There exists a finite Galois extension $L$ of $F$ such that all
endomorphisms of $X$ are defined over $L$. Clearly,
$\Gal(L)=\Gal(F_a/L)$ is an open normal subgroup of finite index
in $\Gal(F)$ and
$$H'=\rho_{2,X}(\Gal(L)) \subset \Aut_{\Z_2}(T_2(X))\subset \Aut_{\Q_2}(V_2(X)))$$
is an open normal subgroup of finite index in $H$.
 We write $\End^0(X)$ for
the $\Q$-algebra $\End(X)\otimes\Q$ of endomorphisms of $X$.

Recall \cite{MumfordAV} that the natural map
$\End^0(X)\otimes_{\Q}\Q_{2} \to \End_{\Q_{2}}V_{2}(X)$ is an
embedding, whose image lies in the centralizer
$$\End_{\Gal(L)}V_{2}(X)=\End_{H'}V_{2}(X).$$
Since the $H$-module $V_2(X)$ is very simple and $H'$ is normal in
$H$, we conclude, thanks to Remark \ref{image}(iv) that either the
$H'$-module $V_2(X)$ is absolutely simple or $H'$ consists of
scalars. In the former case $\End_{H'}V_{2}(X)=\Q_2$ and therefore
$$\End^0(X)\otimes_{\Q}\Q_{2}\cong \Q_2.$$ This implies easily
that $\End(X)=\Z$ and we are done. So, let us assume that $H'$
consists of {\sl scalars}. We are going to arrive to a
contradiction that proves the Lemma.

Since $H' \subset H$ consists of symplectic automorphisms of
$V_2(X)$,
 either $H=\{1\}$ or $H=\{\pm 1\}$.
 In particular, $H'$ is always finite.
Since $H'$ is a subgroup of finite index in $H$, the group $H$ is
finite. In particular, the kernel of the reduction map modulo $2$
$$\Aut_{\Z_2}T_2(X) \supset H \to G \subset \Aut(X_2)$$
consists of periodic elements and, thanks to Minkowski-Serre Lemma
\cite{SZ}, $Z:=\ker(H \to G)$ has exponent $1$ or $2$. This
implies that  $Z$ is commutative.

 Since $Z$ is normal in $H$ and the $H$-module $V_2(X)$ is very
 simple,
either $Z$ consists of scalars  or  the $Z$-module $V_2(X)$ is
absolutely simple, thanks to Remark \ref{image}(iv). On the other
hand, since $Z$ is commutative and $\dim_{\Q_2}(V_2(X))=2g>1$, the
$Z$-module $V_2(X)$ is {\sl not} absolutely simple. Hence $Z$
consists of scalars. Since $Z \subset H$ consists of symplectic
automorphisms of $V_2(X)$,
 either $Z=\{1\}$ or $Z=\{\pm 1\}$.
In other words, either $H=G$ or $H$ is a central double cover of $G$.
As we have already seen, $V_2(X)$ is an absolutely irreducible symplectic $2g$-dimensional representation of $H$.
However, if $G=\Sz(8)$ then, according to the Atlas ~\cite[p. 28]{Atlas},
 $H$ does not have an absolutely irreducible {\sl symplectic} $64$-dimensional representation in characteristic zero.
  Therefore $G$ is not isomorphic to $\Sz(8)$. Also, if
  $G=\L_4(3)$ then, according to the Atlas ~\cite[pp.
  68--69]{Atlas}, $H$ does not have an absolutely irreducible  $38$-dimensional representation in characteristic zero.
Therefore $G$ is not isomorphic to $\L_4(3)$
    and, by \ref{supernot}(b),  the Schur multiplier of $G$ has odd order.
  This implies that $H$ is a trivial central
extension of $G$, i.e., either $Z=\{1\}$ and $H=G$ or $Z=\{\pm
1\}$ and $H=G\times \{\pm 1\}$. In both cases one may view
$V_2(X)$ as absolutely simple symplectic $\Q_2[G]$-module and also
as a lifting of the $\F_2[G]$-module $X_2$ to characteristic zero.
The condition \ref{supernot}(b) means that the representation of
$G$ in $V_2(X)$ satisfies the condition of classical
Brauer-Nesbitt theorem (~\cite[Theorem 1]{Brauer}; see also
~\cite[\S 62]{D}, ~\cite[p. 249]{Hump}) and therefore has {\sl
defect} $0$ ~\cite[\S 16.4]{Serre}. This implies that its {\sl
reduction}, i.e., the $\F_2[G]$-module $X_2$, is {\sl projective}
and all  liftings  of $X_2$ to characteristic zero must be
isomorphic ~\cite[\S 14.4]{Serre}. Since, by \ref{supernot}(d),
there exists an orthogonal lifting of $X_2$, we conclude that the
$G$-module $V_2(X)$ must be also orthogonal. But $V_2(X)$ is
symplectic absolutely simple and therefore cannot be orthogonal.
This gives us the desired contradiction.

\section{Complements to \cite{ZarhinMMJ}}
\begin{thm}
\label{L34} Let $K$ be a field with $\fchar(K) \ne 2$,
 $K_a$ its algebraic closure,
$f(x) \in K[x]$ an irreducible separable polynomial of  degree
$n$. Let us assume  that $n$ and the Galois group $\Gal(f)$ of $f$
enjoy  the following properties:
\begin{enumerate}
\item[(i)]
There exist a positive integer $m>2$ and an odd power prime $q$
such that $n=\frac{q^m-1}{q-1}$;
\item[(ii)]
 $\Gal(f)$ contains a subgroup isomorphic  to $\L_m(q)$.
\end{enumerate}
Let $C_f$ be the hyperelliptic curve $y^2=f(x)$. Let $J(C_f)$ be
its jacobian, $\End(J(C_f))$ the ring of $K_a$-endomorphisms of
$J(C_f)$. Then $\End(J(C_f))=\Z$.
\end{thm}

\begin{proof}
When $(q,m)\ne (3,4)$, the assertion of Theorem \ref{L34} is
proven in \cite{ZarhinMMJ}. (The case $(q,m)= (3,4)$ was treated
in \cite{ZarhinP} under an additional assumption that
$\fchar(K)=0$.) So, in the course of the proof we may assume that
$(q,m)= (3,4)$ and therefore $n=40, g=19$. Replacing $K$ by a
suitable separable extension, we may assume that
$G:=\Gal(f)=\L_4(3)$. Taking into account that $\L_4(3)$ is simple
non-abelian and replacing $K$ by its abelian extension obtained by
adjoining to $K$ all $2$-power roots of unity, we may also assume
that $K$ contains all $2$-power roots of unity.

Let $\RR$ be the set of roots of $f$. It follows from Theorem
\ref{Lnq} that the $G=\Gal(f)$-module $Q_{\RR}$ is very simple. It
follows from Corollary 5.3 of \cite{ZarhinTexel} that either
$\End(J(C_f))=\Z$ or  $\fchar(K)>2$ and $J(C_f)$ is a
supersingular abelian variety. Now the result follows from Lemma
\ref{supernot}.
\end{proof}

\begin{ex}
Suppose $p$ is an odd prime, $q>1$ is a power of $p$, $m>2$ is an
even integer. Let us put $n=(q^m-1)/(q-1)$.  Suppose $k$ is an
algebraically closed field of characteristic $p$ and $K=k(z)$ is
the field of rational functions. The Galois group of $x^n+zx+1$
over $K$ is $\L_m(q)$ and the Galois group of $x^n+x+z$ over $K$
is $\PGL_m(\F_q)$ \cite[p. 1643]{Ab}. Therefore the jacobians of
the hyperelliptic curves $y^2=x^m+zx+1$ and $y^2=x^m+x+z$ have no
nontrivial endomorphisms over an algebraic closure $K_a$ of $K$.
It follows from Proposition 4.5 and Example 4.2(iv) of
\cite{ZarhinSh} that these jacobians are {\sl not} isogenous over
$K_a$ if $m$ and $q-1$ are {\sl not} relatively prime .
\end{ex}

\end{document}